\documentclass[10pt, twocolumn, conference]{IEEEtran}

\usepackage{amsmath, mathrsfs,amsfonts}
\usepackage{amssymb}
\usepackage{graphicx}
\usepackage{cite}
\usepackage{algorithm2e}
\usepackage{array}
\usepackage[tight]{subfigure}
\usepackage[top=0.74in, bottom=0.74in, left=0.55in, right=0.55in]{geometry}

\usepackage{enumerate}

\IEEEoverridecommandlockouts
\begin{document}
\title{Modeling the Impact of Communication Loss on the Power Grid under Emergency Control \thanks{This work was supported by DTRA grant HDTRA1-13-1-0021.}}

\author{\IEEEauthorblockN{Marzieh Parandehgheibi \IEEEauthorrefmark{1}, Konstantin Turitsyn \IEEEauthorrefmark{2} and Eytan Modiano \IEEEauthorrefmark{1}}
\IEEEauthorblockA{\IEEEauthorrefmark{1}Laboratory for Information and Decision Systems, Massachusetts Institute of Technology, Cambridge, MA, USA\\
\IEEEauthorrefmark{2}Department of Mechanical Engineering, Massachusetts Institute of Technology, Cambridge, MA, USA\\}}

\maketitle 

\begin{abstract}
We study the interaction between the power grid and the communication network used for its control. We design a centralized emergency control scheme under both full and partial communication support, to improve the performance of the power grid. We use our emergency control scheme to model the impact of communication loss on the grid. We show that unlike previous models used in the literature, the loss of communication does not necessarily lead to the failure of the correspondent power nodes; i.e. the ``point-wise" failure model is not appropriate. In addition, we show that the impact of communication loss is a function of several parameters such as the size and structure of the power and communication failure, as well as the operating mode of power nodes disconnected from the communication network. Our model can be used to design the dependency between the power grid and the communication network used for its control, so as to maximize the benefit in terms of intelligent control, while 
minimizing the risks due to loss of communication.

\vspace{-1mm}
\end{abstract}

\IEEEpeerreviewmaketitle

\section{Introduction}\label{Intro_sec}

In August 2003, 50 million people lost power due to the large blackout in North East America. According to the aftermath reports \cite{FinalReport2004:Online,NERC2003:Online,NY-ISO2004:Online}, initially, three power lines in Ohio were disconnected due to inadequate tree trimming. These failures caused tripping of overloaded lines and shedding of generators experiencing under-frequency; however, these changes were not monitored by grid operators. Consequently, the large imbalance in the power grid led to a catastrophic cascade of failures in other regions of United States and Canada. Similarly, in September 2011, 2.7 million people lost power in Arizona, Southern California, and Baja California, Mexico. This major blackout started by tripping of a single power line; however, since the grid was not operating in the secure N-1 state, it led to a rapid cascade of failures in power lines and generators \cite{CaliforniaReport2012:Online}. Studies \cite{FinalReport2004:Online,CaliforniaReport2012:Online,andersson2005causes} show that the main reasons behind both blackouts were a lack of situational awareness and a lack of coordination between grid operators in neighbor regions.

During normal operation, primary and secondary frequency controls are responsible for stabilizing the grid. In particular, primary frequency control is a local controller which reacts to local changes in frequency and adjusts the generation to keep the frequency within an acceptable range. The secondary controller is responsible for setting the frequency back to its nominal value (e.g. $60 Hz$ in US) where it uses the generator's reserves to balance the power. However, during large failures the normal operation controllers cannot stabilize the grid. Therefore, the future smart grid should be equipped with a Communication and Control Network (CCN) that allows rapid monitoring of the power grid and provides intelligent centralized control actions that can mitigate cascade of failures. The centralized control actions to stabilize the grid during catastrophic failures are referred to as ``Emergency Control''.

There are several studies proposing Emergency control schemes for changing the power generation as well as load shedding so that the grid can be stabilized before any cascading failures occur \cite{chen2005cascading,bienstock2011optimal,koch2010mitigation,pahwa2013load}. Although, using this extra information/control improves the performance of the grid, it creates a dependency between the power grid and the communication and control network.

During normal operations, loss of communication is unlikely to lead to significant power failures as local controllers can stabilize the grid. However, when the grid is under stress, lack of situational awareness and control can lead to catastrophic failures. Such events may result from a natural disaster that affects both the communication network and the power grid or the failure of communication components due to loss of power coming from the grid. Therefore, it is very important to study the impact of communication loss on the power grid's stability during a large disturbance.
 
\begin{figure}[h]
\centering
\subfigure[No Control- Interdependent networks are more vulnerable.]
{\label{VulnerableInterdep}\includegraphics[scale=0.04]{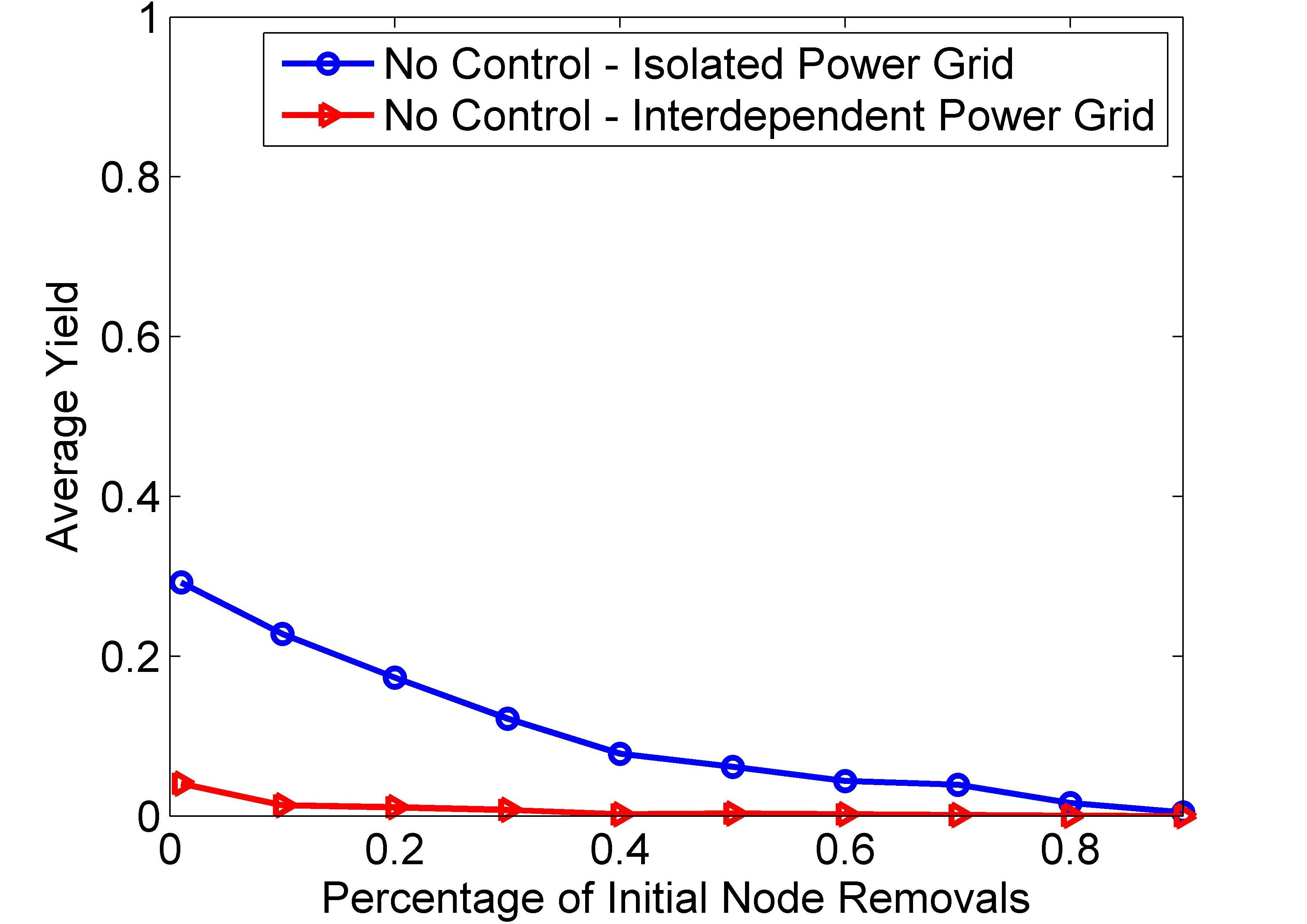}}                
\subfigure[Control - Interdependent networks are more reliable.]
{\label{ImprovedInterdep}\includegraphics[scale=0.04]{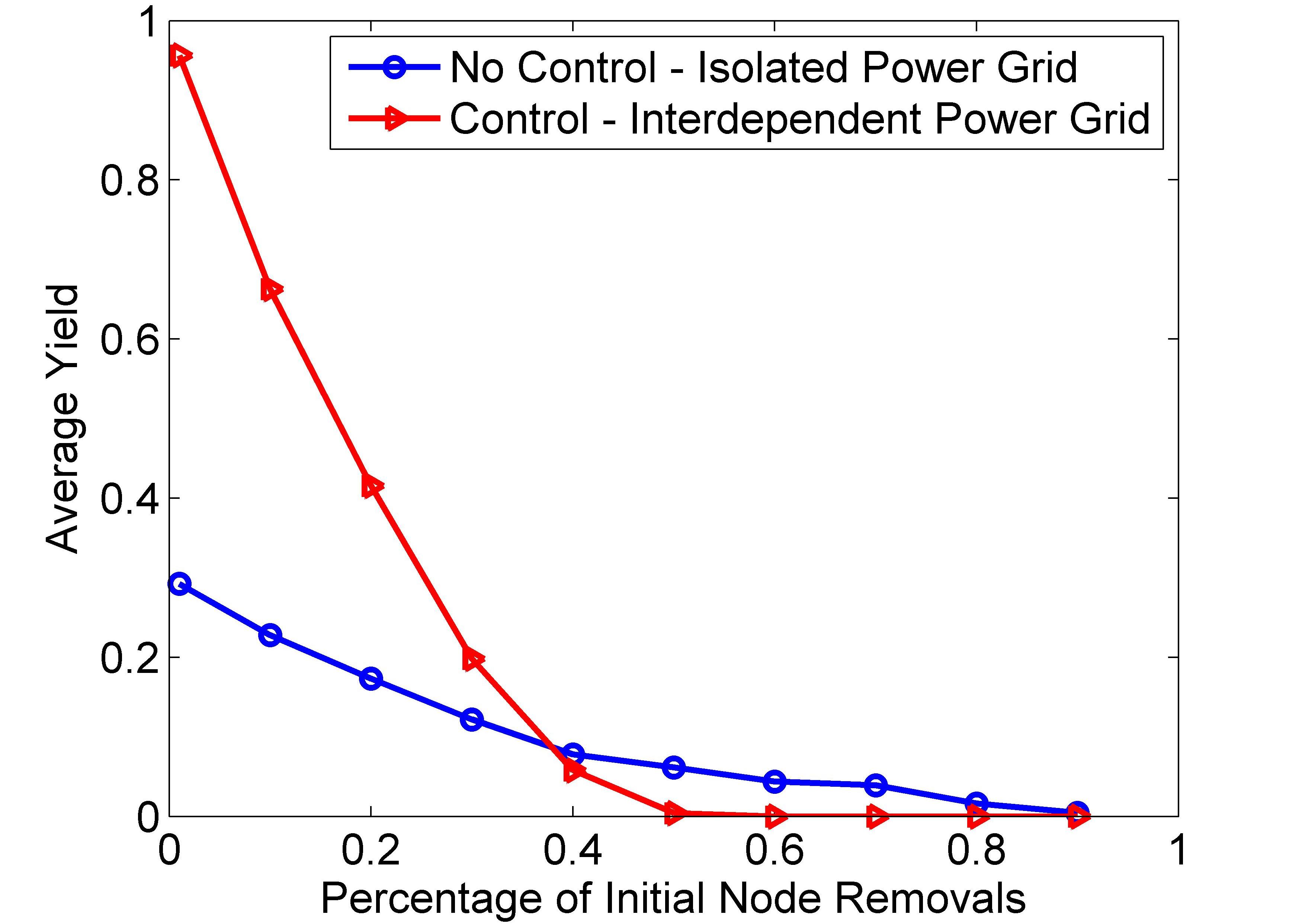}}
\vspace{-1mm}
\caption{\small{Comparing yield of interdependent and isolated power grids.}}
\label{Comparison_Interdep_Isolated}
\end{figure}

The impact of communication on the power grid's performance and vice versa was recently studied using an abstract form of interdependency. Buldyrev \textit{et al.} in \cite{Buldyrev2010} showed that if there exists a one-to-one interdependency between the nodes of the grid and the communication network, interdependent networks are more vulnerable to failures than isolated networks. In their ``point-wise" interdependency model, a power node fails if it loses its connection to the communication network, and a communication node fails if it loses its connection to the power grid. Similar results were obtained in \cite{gao2012networks,parandehgheibi2013robustness}. Recently, \cite{parandehgheibi2014mitigating} showed that it is critical to use the power flow equations to model the power grid, and interdependency could benefit the power grid if the communication network is used for mitigating the cascade of failures (See Figure \ref{Comparison_Interdep_Isolated}). As can be seen, when the yield in interdependent power grid and communication without control is lower than the isolated power grid (Figure \ref{VulnerableInterdep}). However, when intelligent control is applied to the interdependent network, the yield is higher than the isolated power grid (Figure \ref{ImprovedInterdep}).
\vspace{-1mm}

In this paper, we carefully model the function of emergency control of the grid using the communication and control network. Using this model, we show that the loss of any communication network component may lead to the loss of situational awareness and control, and impact the performance of the grid. In particular, we show that the point-wise failure model is not suitable for modeling interdependent power-communication networks. Moreover, we show that the impact of communication loss on the power grid is a function of several parameters such as size and structure of the communication and power failures.

The rest of this paper is organized as follows. In Section \ref{Model_Sec} we describe the power-communication dependency. In Section \ref{Emergency_Sec}, we formulate the emergency control problem with full and partial communication. Finally, in Section \ref{Simulation_Sec}, we provide simulation results and conclude in Section \ref{Conclusion_Sec}.

\section{Model}\label{Model_Sec}

\subsection{PowerGrid-Communication Dependency}
\begin{figure}[h]
\centering
\includegraphics[scale=0.35]{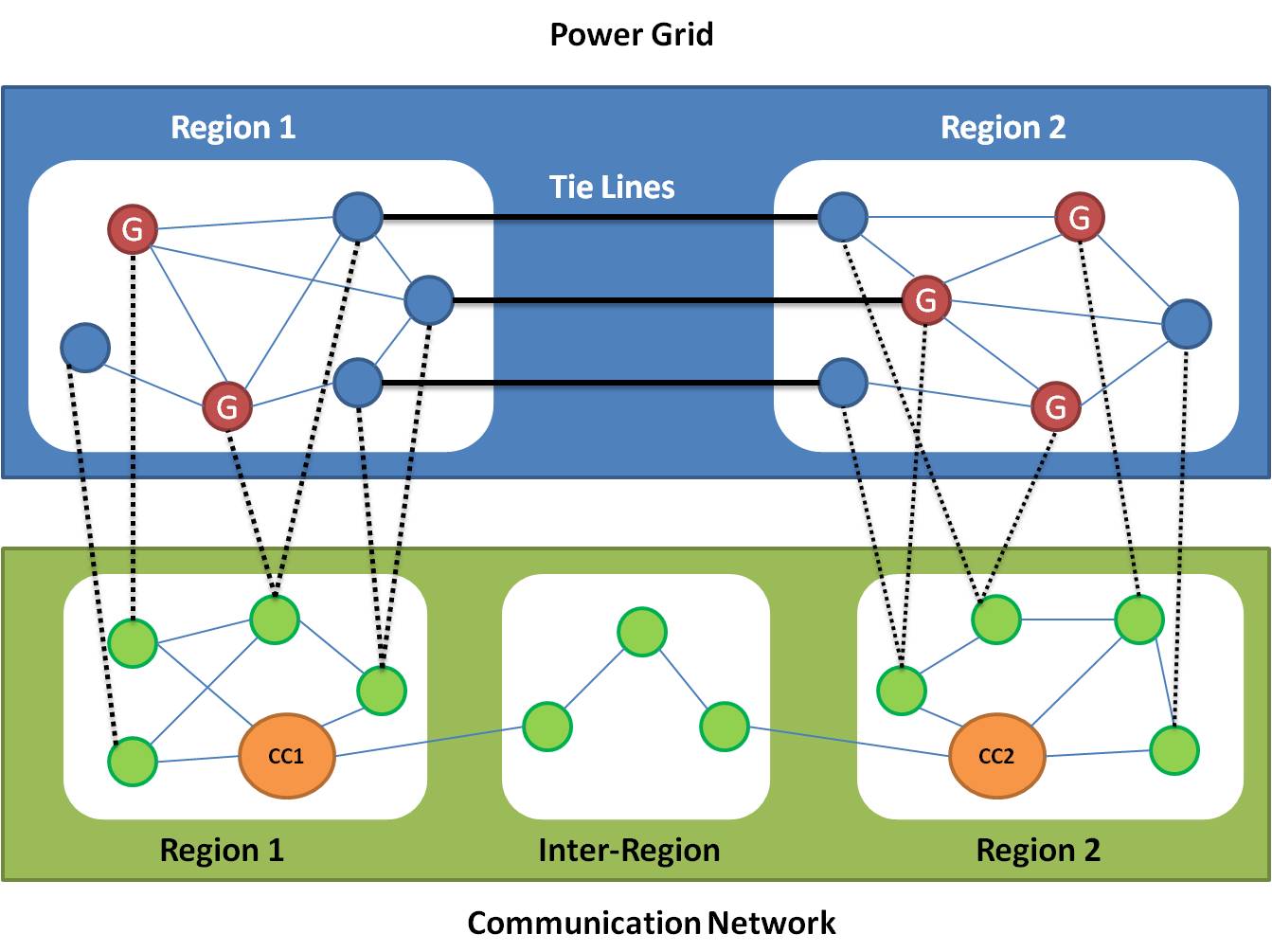}
\caption{Future Power Grid equipped with communication networks and control centers for online monitoring and control; dotted lines indicate dependency of power nodes on communication nodes. Power nodes that have lost their connection to the communication network become ``uncontrollable".}
\label{PowerCommStructure}
\end{figure}

Figure \ref{PowerCommStructure} shows an abstract model of the future control and communication network for an interconnected grid. In this model, each region is supported by a dedicated communication network (intra-region communication network) that monitors all parts of the grid, and sends the information to the control center. The control decisions made by the control center are sent back to the grid via the same communication network. Moreover, the grid regions are connected to each other via several tie lines that allows transmission of power from one region to the other. Therefore, a failure in one region could cascade to the other regions. In order to avoid such failures, the control centers are connected to each other with an inter-region communication network that allows them to share state information and regional control decisions.

It can be seen from Figure \ref{PowerCommStructure} that communication failures could occur either between the regions or inside the regions. \textbf{Inter-region} failures can degrade or disconnect the communication between the control centers in different regions. In contrast, \textbf{Intra-region} failures can cause some of the power nodes to be disconnected from the communication network and unable to send information to or receive the control decisions from the control center.

\begin{figure}
\centering
\includegraphics[scale=0.35]{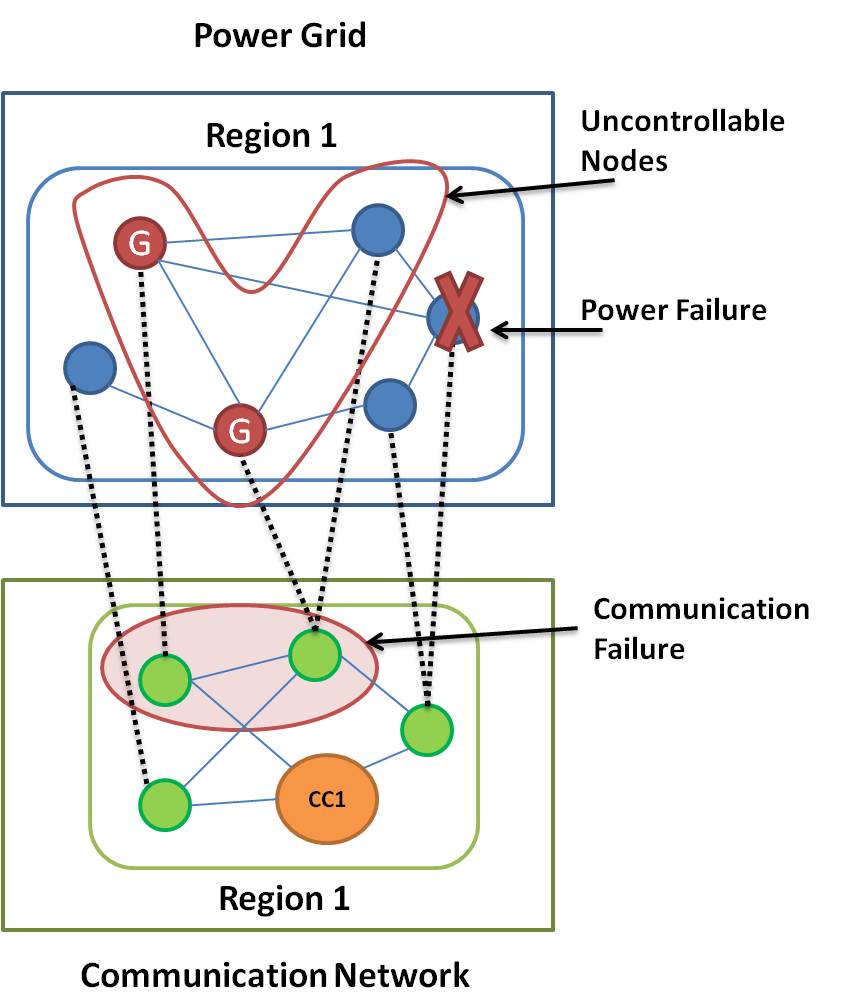}
\caption{Dotted lines indicate dependency of power nodes on communication nodes. Power nodes that have lost their connection to the communication network become ``uncontrollable".}
\vspace{-5mm}
\label{Inter_regionFailure}
\end{figure}

In this paper, we study the case of intra-region failures. In particular, we focus on the case where we are able to observe the power failures; i.e. we have communication to the failed power nodes (See Figure \ref{Inter_regionFailure}).

\subsection{Power Grid}
The power grid consists of nodes and power lines where power nodes have three types: Generators $(G)$ that generate power; Loads $(L)$ that consume power and Buses $(B)$ that allow the transmission of power through them, but neither generate nor consume power. In our model, each generator or load is connected to the rest of the grid via a single bus and each bus can be connected to any number of buses. Figure \ref{PowGrid_CommLoss} shows an example of our system. 

In this paper, we study the steady-state behavior of the grid. We also assume that power lines are loss-less; i.e. can be modeled only by their reactance $X$. Thus, the flow in a power line is described by a DC model as follows: $f_{ij}X_{ij}=\Delta \theta_{ij}$ where $f_{ij}$ is the amount of power in line $(i,j)$, $X_{ij}$ is the reactance of line $(i,j)$ and $\Delta \theta_{ij}$ represents the difference in the voltage phase of nodes $i$ and $j$. 

The generators are modeled as synchronous generators with the following swing equation:
\begin{equation}
\resizebox{0.9\hsize}{!} {$M_i \dot{\omega_i} = PG_i - D_i(\omega_i - \omega_s) - (\sum_{j \in E} f_{ij} - \sum_{j \in E} f_{ji})$}
\end{equation}

where $PG_i$ denotes the mechanical power, and $D_i$ and $M_i$ are generator's damping and inertia, respectively. Moreover, the difference in the flow passing through the generator ($\sum_{j \in E} f_{ij} - \sum_{j \in E} f_{ji}$) denotes the amount of generated electrical power. Finally, $\omega_i$ is the frequency of node $i$ and $\omega_s$ is the synchronous frequency [see \cite{chakrabortty2013introduction} for more details]. Considering the steady-state behavior ($\dot{\omega_i}=0$), the equation for synchronous generators will reduce to:
\begin{equation}
PG_i - D_i(\omega_i - \omega_s) = \sum_{j \in E} f_{ij} - \sum_{j \in E} f_{ji} \label{Sync_eqn}
\end{equation}

The damping coefficient $D_i$ is often $0.02$ per unit or less (See \cite{glover2011power}-pp. 657-663 for more details). This model describes the reaction of generator to the changes in power; if generation is greater than the load, the frequency increases, and if the generation is smaller than the load, the frequency decreases. In order to avoid sharp changes in frequency, we also model the local frequency-droop control (also called primary frequency control) in the grid that is responsible for modifying the amount of power generation based on the changes in frequency; i.e. increases the power generation as frequency drops and decreases the power generation as frequency rises. The droop-control can be written as follows.

\begin{equation}
\Delta PG_i = -\frac{1}{R_i}(\omega_i - \omega_s) \label{DroopCntl_eqn}
\end{equation}

where $R_i$ is the regulation constant. In a per unit system, the standard value of regulation constant is $0.05$ per unit. Thus, in our model of a $60 Hz$ system, $R_i=\frac{60*0.05}{PG_i^{initial}}$ where $PG_i^{initial}$ is the initial mechanical power of generator $i$ before any disturbance (\cite{glover2011power}-pp. 657-663). The combination of $PG_i$ in equation \ref{Sync_eqn} and $\Delta PG_i$ in equation \ref{DroopCntl_eqn} models the response of generator $i$ to the changes in power using the local controller:
\begin{equation}
\begin{aligned}
& PG_i -(D_i+\frac{1}{R_i})(\omega_i - \omega_s) & = \sum_{j \in E} f_{ij} - \sum_{j \in E} f_{ji} \label{Final_Sync_eqn}
\end{aligned}
\end{equation}

For the rest of this paper, we define $\alpha_i=D_i+\frac{1}{R_i}$.

Finally, we consider protection relays located at every power node. The control role of these relays is explained in section \ref{control_actions_sec}.

\subsection{Communication Network}
We consider a zero-delay communication network that supports every power node; i.e. collects synchronous information from every power node, sends it to the control center and sends back the control decision from the control center to all power nodes. The set of collected information by the communication network includes: (i) Magnitude of voltage at node $k$ ($V_k(t)$); (ii) Phase of voltage at node $k$ ($\theta_k(t)$); (iii) Frequency of node $k$ ($\omega_k(t)$); (iv) Flow in power line $(k,j)$ ($f_{kj}(t)$); (v) on/off State of element $j$ ($S_j(t)$).

Note that since we use the DC model in the paper, the magnitude of voltage $V_k(t)$ is constant and equal to 1 for all nodes. Moreover, the ``off" state of an element means that it has failed.

\subsection{Control Actions}\label{control_actions_sec}
Next, we describe both types of local and centralized control actions needed for operation of the power grid.

\textbf{Local Control Actions}: All power nodes are equipped with local controllers that do not require connection to the communication network and their actions include: (i) \textit{droop control at generators}: Droop control can increase or decrease the amount of generation based on the changes in the frequency as described by equations \ref{DroopCntl_eqn} and \ref{Final_Sync_eqn}; (ii) \textit{Over-frequency generator tripping (protection)}: the protection relays will trip the generator if the frequency exceeds the maximum threshold $\omega_{max}$ due to excess generation in the system; (iii) \textit{Under-frequency load shedding (protection)}: the protection relays will shed the load if the frequency drops below the minimum threshold $\omega_{min}$ due to excess load in the system; (iv) \textit{Overloaded line tripping (protection)}: the protection relays will trip the line if the power in that line exceeds the capacity.

\textbf{Central Control Actions}: All power nodes are equipped with sensors/actuators connected to the control center via the communication network. The central control actions include: (i) \textit{Ramping down generators}: if power generation is greater than consumption, the controller decreases the generation to keep the frequency within the acceptable range; (ii) \textit{Intelligent load shedding}: if the power generation is lower than consumption, the controller sheds some load to keep the frequency within the acceptable range; (iii) \textit{Intelligent line tripping}: can be used for changing the topology of the grid or islanding some areas of the power grid. Since in our model every power node is equipped with an actuator, a power line can be tripped by either of its end-nodes.

\section{Emergency Control}\label{Emergency_Sec}
In this Section, we design the optimal emergency control used for mitigation of failures in the power grid in the presence of a fully or partially operational communication network. We can then use these control policies for evaluating the performance of the grid under different communication failure scenarios.

\subsection{Full Communication}\label{Full_Emergency}
In this case, the power grid is fully supported by a communication network; thus, every generator, load and line can be centrally controlled. The optimal emergency control is the set of central and local control actions that maximizes the served load while keeping the power balanced, maintaining the frequency within an acceptable range and keeping the flows in the power lines within their capacity. We use the control model with full communication as a basis for modeling the control with partial communication.

In eq. (\ref{Full_EC_Formulation}), we formulate the problem of optimal control with the objective of maximizing served load while stabilizing the grid. Let, $V_G$, $V_L$ and $V_B$ denote the set of generators, loads and bus nodes after the initial power failure, respectively. Moreover, let $\omega_i^{min}$ and $\omega_i^{max}$ denoted the minimum and maximum frequency thresholds. In addition, $E$ denotes the set of power lines after the initial failure, and $X_{ij}$ represents the reactance of line $(i,j)$. Finally, let $M$ be a large constant.

The variables $f_{ij}$ and $\Delta \theta_{ij}$ denote the amount of flow in line $(i,j)$ as well as the phase difference of voltages at nodes $i$ and $j$. In addition, variables $PG_i$, $PL_i$ and $\omega_i$ denote the amount of generation, load and frequency at node $i$. Note that $PG_i$ and  $PL_i$ take positive and negative values, respectively. Finally, $z_{ij}$ is a binary variable associated to line $(i,j)$ that takes value of $1$ if line $(i,j)$ is connected and $0$ if that line is tripped (modeled in constraint (\ref{Full_Trip_Line})).

As mentioned previously, the objective is to serve the maximum load while stabilizing the grid. Constraints (\ref{Full_Gen},\ref{Full_Load},\ref{Full_Bus}) ensure that the flow conservation is satisfied in generators, loads and buses where generation and load values $PG_i$ and $PL_i$ can change by central controller and term $\alpha_i (\omega_i-\omega_s)$ models the local droop controller at the generators. Constraint (\ref{Full_flow}) models the DC power flow in line $(i,j)$ if it remain connected; i.e. $z_{i,j}=1$; note that there is no relation between the phases at nodes $i$ and $j$ if the line is tripped ($z_{i,j}=0$). Constraint (\ref{Full_capacity}) guarantees that the flow in line $(i,j)$ is within the capacity if the line is not tripped, and constraint (\ref{Full_frequency}) ensures that all nodes in a connected area have the same frequency. Finally, constraints (\ref{Full_freq_cap}, \ref{Full_Gen_cap}, \ref{Full_Load_cap}) guarantee that the values of frequency, generation and load are maintained within the acceptable range.

\vspace{-1mm}
\begin{subequations}
\small{
\begin{alignat}{3}
\min & \sum_{i \in V_L}PL_i  \label{Full_obj}\\
& \sum_{j \in E} f_{ij} - \sum_{j \in E} f_{ji} = PG_i - \alpha_i (\omega_i-\omega_s) && \forall i \in V_G \label{Full_Gen}\\
& \sum_{j \in E} f_{ij} - \sum_{j \in E} f_{ji} = PL_i && \forall i \in V_L \label{Full_Load}\\
& \sum_{j \in E} f_{ij} - \sum_{j \in E} f_{ji} = 0 && \forall i \in V_B \label{Full_Bus}\\
-& M(1-z_{ij}) \leq X_{ij} f_{ij} -\Delta \theta_{ij} \leq M(1-z_{ij}) && \forall (i,j) \in E \label{Full_flow}\\
-& z_{ij}f_{ij}^{max} \leq f_{ij} \leq z_{ij}f_{ij}^{max} && \forall (i,j) \in E \label{Full_capacity}\\
-& M(1-z_{ij}) \leq \omega_i - \omega_j \leq M(1-z_{ij}) && \forall (i,j) \in E \label{Full_frequency}\\
& \omega_i^{min} \leq \omega_i \leq \omega_i^{max} && \forall i \in V_G \label{Full_freq_cap}\\
& PG_i^{min} \leq PG_i - \alpha_i (\omega_i-\omega_s) \leq PG_i^{max} && \forall i \in V_G \label{Full_Gen_cap}\\
& PL_i^{max} \leq PL_i \leq 0 && \forall i \in V_L \label{Full_Load_cap}\\
& z_{ij} \in \{0,1\} && \forall (i,j) \in E \label{Full_Trip_Line}
\end{alignat}
\vspace{-5mm}
}\label{Full_EC_Formulation}
\end{subequations}
\vspace{-1mm}

\subsection{Partial Communication}\label{Partial_Emergency}

Next, we consider the case that in addition to the power failure, a part of the communication network fails as well. Thus, parts of the grid lose their connection to the communication network and control center. Our objective is to design an emergency control policy that maximizes the served load and stabilizes the grid using the controllable nodes.

\begin{figure}[h]
\centering
\includegraphics[scale=0.50]{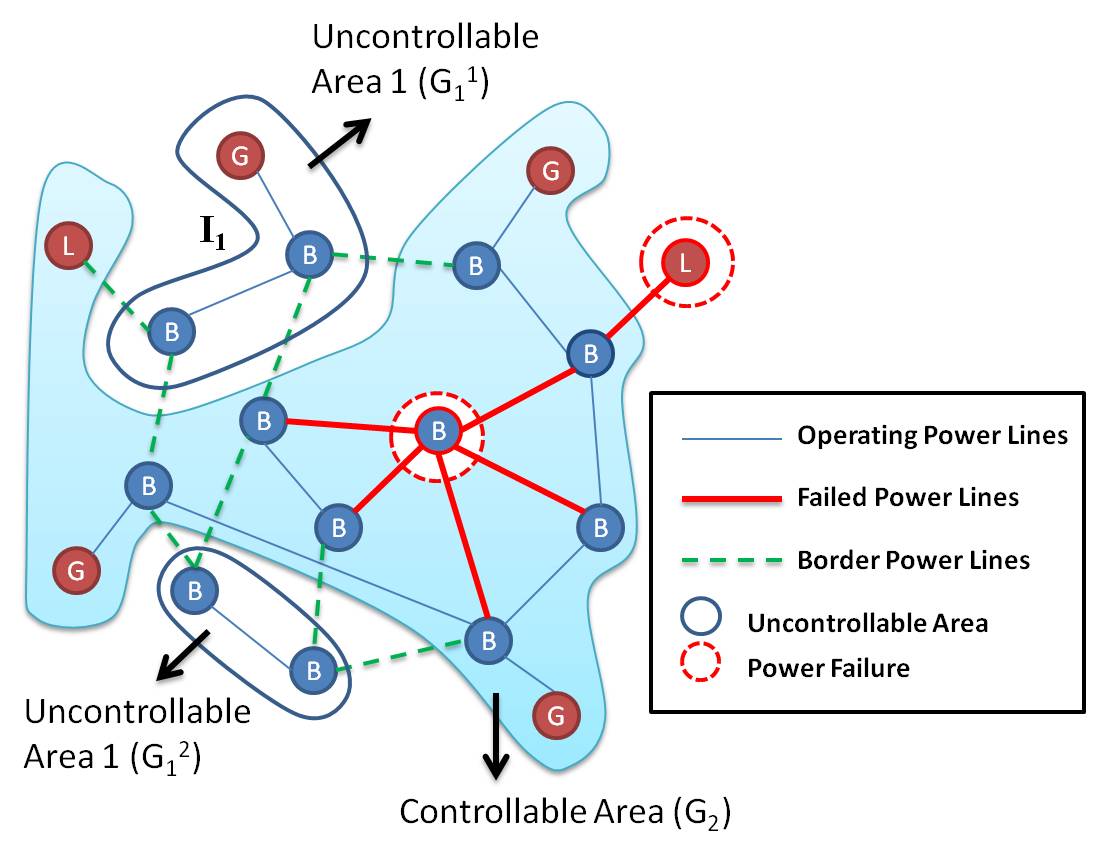}
\caption{Power Grid Model after Communication Loss and Power Failure}
\label{PowGrid_CommLoss}
\end{figure}

Figure \ref{PowGrid_CommLoss} shows an example of such power grid with power failures and controllable/uncontrollable areas. In the following, we define each area mathematically, and explain the set of control actions available in each.

The components that have not initially failed in the power grid can be divided into 3 areas described as follows:

\textit{Uncontrollable Area:} Let $G_1=\{V_1,E_1\}$ denote the uncontrollable subset of power grid, where $V_1$ denotes the set of power nodes that are disconnected from the communication network, and $E_1$ denotes the set of power lines that connect any pair of nodes in $V_1$. Without loss of generality, we can divide subgraph $G_1$ into $k$ disjoint subgraphs $G_1^k=\{V_1^{k},E_1^{k}\}$ such that nodes inside $V_1^{k_i}$ are connected and for any $j \neq i$ nodes in $V_1^{k_i}$ and $V_1^{k_j}$ are not connected. Note that depending on the structure of the network, $k$ could be any number equal or greater than 1. The only possible centralized control in these areas is islanding of the entire area by tripping border lines. All of the local controllers described in Section \ref{control_actions_sec}; i.e. droop control and protection schemes, are available. In addition, we have a local controller that switches the nodes in the uncontrollable areas to a pre-defined mode of operation as soon as they are disconnected from the communication network. In this paper, we consider two possible modes: $P_{init}$ which corresponds to keeping all nodes operating at their last state, and $P_{zero}$ which corresponds to tripping all nodes.

\textit{Controllable Area:} Let $G_2=\{V_2,E_2\}$ denote the controllable subset of the power grid, where $V_2$ is the set of power nodes that are connected to the communication network ($V_2 \cap V_1 =\emptyset$ and $V_2 \cup V_1 =V$  where $V$ is the set of all power nodes), and $E_2$ is the set of power lines connecting any pair of nodes in $V_2$. All central and local controllers described in Section \ref{control_actions_sec} are available in the controllable areas.

\textit{Border Lines:} Let $E_{12}^k$ denote the set of power lines that connect the uncontrollable nodes in $V_1^k$ to the controllable nodes in $V_2$. Border lines can be tripped centrally using the relay in its controllable node. All local controllers at its end-nodes are available.

In eq. (\ref{Partial_EC_Formulation}), we formulate the optimal emergency control problem subject to loss of communication. Given the set of control actions, the behavior of the controllable areas can be modeled as described in previous section with Full Control; i.e., eq. (\ref{controllable_const}). However, the nodes in the uncontrollable areas cannot employ centralized control and must rely on localized control and/or islanding as described next. The following ILP shows the general description of optimal emergency control.
\vspace{-3mm}

\begin{subequations}
\small{
\begin{alignat}{3}
\min & \sum_{i \in V_L}PL_i  \label{Partial_obj}\\
\mbox{s.t.} & \mbox{  Constraints} (\ref{Full_Gen}-\ref{Full_Trip_Line}) 			&& \quad \forall k \in K, i \in V_2, (i,j) \in E_2 \cup E_{12}^k \label{controllable_const}\\
& \mbox{  Constraints} (\ref{island_tripping_relation}-\ref{stableIsland}) 			&& \quad \forall k\in K, i \in V_{1}^k,  (i,j) \in E_1^k
\end{alignat}
}
\label{Partial_EC_Formulation}
\end{subequations}
 
\vspace{-3mm}
In the following, we describe constraints (\ref{island_tripping_relation}-\ref{frequency_const_uncontrollable}) related to the control actions in the uncontrollable areas. Let $\sum_{(i,j) \in E_{12}^k} z_{ij}=0$ if the uncontrollable area $G_1^k$ is isolated from the rest of grid; i.e. all the border lines are tripped, and $\sum_{(i,j) \in E_{12}^k} z_{ij}>0$ if it is connected. Moreover, let $I^k$ be a binary variable (modeled in constraint (\ref{stableIsland})) associated with the uncontrollable area $G_1^k$ where $I^k=0$ if the isolated uncontrollable area cannot be stabilized just by the local droop controller, and $I^k=1$ if that area is stabilized (i.e. power is balanced just by using droop control and frequency and power flows are within the acceptable range). 

Constraints (\ref{island_tripping_relation}-\ref{frequency_const_uncontrollable}) model the control decision in the uncontrollable areas. Note that $PG_i^{init}$ and $PL_i^{init}$ denote the values of generators and loads in uncontrollable area after going to mode $P_{init}$; i.e. the last values of generators and loads right before disconnection from communication network. If the mode is set to $P_{zero}$, $PG_i^{init}$ and $PL_i^{init}$ in constraints (\ref{powerGen_unconrollable}) and (\ref{powerLoad_unconrollable}) will be set to zero. Constraint (\ref{island_tripping_relation}) denotes that if the uncontrollable area is unstable ($I^k=0$), it has to be islanded ($\sum_{(i,j) \in E_{12}^k} z_{ij}$). 

When an uncontrollable area can be stabilized just by local controllers; i.e. $I^k=1$, constraints (\ref{island_tripping_relation}-\ref{frequency_const_uncontrollable}) will be active. In particular, constraints (\ref{powerGen_unconrollable}-\ref{powerBus_unconrollable}) model the power balance in the area using only the droop control at generators, and constraint (\ref{OhmLaw}) models the DC power flow in line $(i,j)$. Moreover, constraint (\ref{linecap}) guarantees that flow is within the line capacities, constraint (\ref{equalFreq}) forces all the connected nodes to have the same frequency and constraint (\ref{frequency_const_uncontrollable}) guarantees that the frequency remains within the acceptable range.

\begin{subequations}
\vspace{-3mm}
\small{\small{
\begin{alignat}{5}
& \sum_{(i,j) \in E_{12}^k} z_{ij} \leq M I^k 					&& \forall k \in K &&\label{island_tripping_relation}\\
-& M(1-I^k) \leq \sum_j f_{ij} - \sum_j f_{ji}- 			PG&&_i^{init} + \alpha_i(\omega_i-\omega_s)&& \nonumber \\
	& \quad\quad\quad\quad\quad \leq M(1-I^k) 						&& \forall i \in V_{G1}^k, \forall k \in K && \label{powerGen_unconrollable}\\
-& M(1-I^k) \leq \sum_j f_{ij} - \sum_j f_{ji}-				PL&&_i^{init}\leq M(1-I^k)  && \nonumber\\	
& \quad\quad\quad\quad\quad   													&& \forall i \in V_{L1}^k, \forall k \in K &&\label{powerLoad_unconrollable}\\
-& M(1-I^k) \leq \sum_j f_{ij} - \sum_j f_{ji}\leq 		M(&&1-I^k)  && \nonumber\\
& 																											&& \forall i \in V_{B1}^k, \forall k \in K && \label{powerBus_unconrollable}\\
-& M(1-I^k)\leq X_{ij} f_{ij}-\Delta \theta_{ij} \leq M(1&&-I^k) && \nonumber\\
&  																											&&\forall (i,j) \in E_1^k, \forall k \in K &&  \label{OhmLaw}\\
- & f_{ij}^{max} I^k \leq f_{ij} \leq f_{ij}^{max} I^k  &&\forall (i,j) \in E_1^k, \forall k \in K &&  \label{linecap}\\
-& M(1-I^k) \leq \omega_i - \omega_j \leq  M(1-I^k)			&&\forall (i,j) \in E_1^k, \forall k \in K && \label{equalFreq} \\
& \omega_i^{min} I^k \leq \omega_i \leq \omega_i^{max} I^k  &&	\forall i \in V_{G1}^k , \forall k \in K && \label{frequency_const_uncontrollable}\\
& I^k \in \{0,1\} 																			&&\forall k \in K&& \label{stableIsland}
\end{alignat}
}}
\vspace{-5mm}
\label{PartialControl_Uncontrollable}
\end{subequations}

Note that when $I^k=0$, the uncontrollable area is unstable; i.e. either power cannot be balanced just by droop controller or some lines are overloaded. Thus, the frequency and line protection relays will be activated to trip the generators, shed the load and trip the lines. This causes extra failures which activates more protection relays. The cascade of failures continue until power is balanced, and frequency and power flows are within their acceptable ranges. We model the cascading failures separately for each uncontrollable area $k$ after observing its correspondent $I^k$ value.

\section{Simulation Results}\label{Simulation_Sec}
We analyze the data from the Italian power grid which consists of 310 buses, 113 generator units and 97 loads. The power failures are considered to be generators. We assume that any arbitrary area in the power grid can lose its connection to the communication network, and investigate the impact of these uncontrollable areas on the performance of the grid. The metric that we use for our analysis is ``yield" defined as the ratio of served load \footnote{Served Load is considered to be the sum of power from controllable areas, stable uncontrollable areas and unstable uncontrollable areas that have experienced cascading failures.} to the initial load.

Simulation results show that loss of communication can indeed impact the performance of the power grid and lead to a lower yield. For example, we observed a scenario where loss of 10 power nodes led to 6\% load shedding under full communication and 23\% load shedding under parial communication (See \cite{parandehgheibiinvestigating}). In this section, we would like to find the parameters that have the greatest impact on performance.

\subsection{Effect of size and structure of uncontrollable areas}
We consider different number of uncontrollable nodes with different clusterings. In particular, we define a cluster as a set of connected nodes where a cluster of size 1 means that no two uncontrollable nodes are connected and a cluster of size 10 means that the uncontrollable area can be divided into disjoint subareas, each with 10 connected nodes.

\begin{figure}[h]
\centering
\includegraphics[scale=0.055]{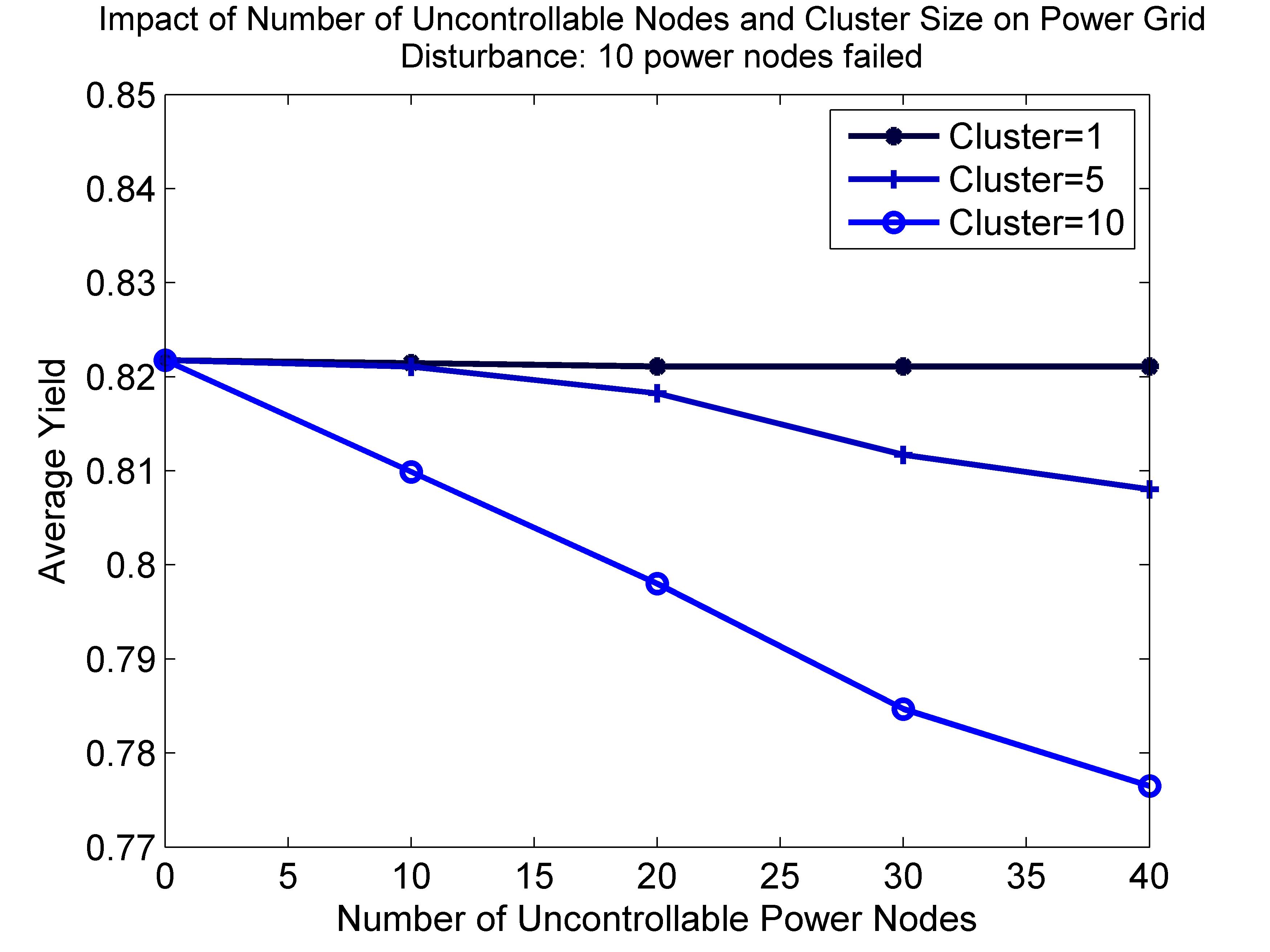}
\caption{Average yield vs. Number of uncontrollable power nodes; $P_{init}$ mode.}
\vspace{-3mm}
\label{figSizeComm}
\end{figure}

\begin{figure}[h]
\centering
\includegraphics[scale=0.055]{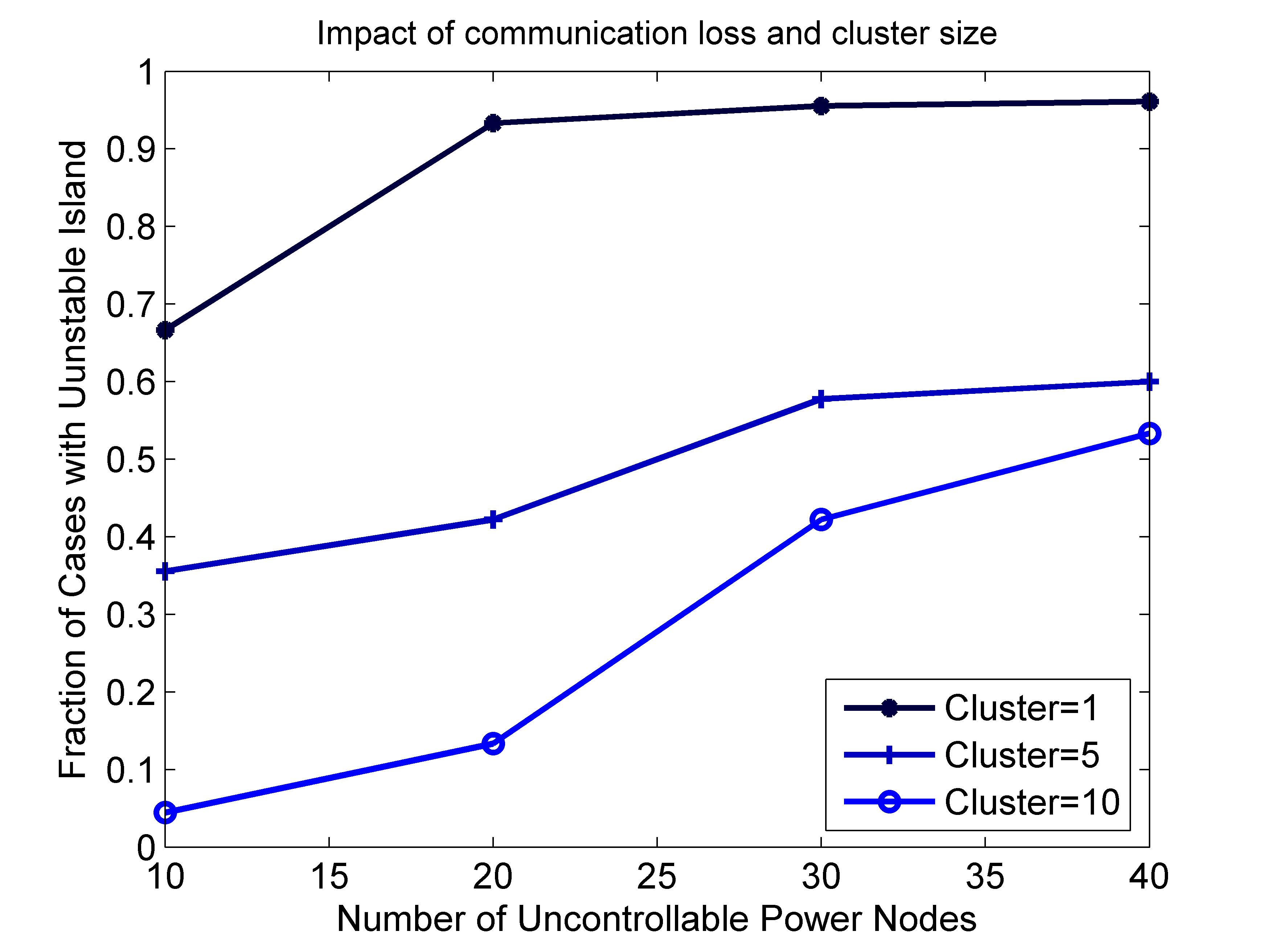}
\caption{Fraction of cases with unstable islands vs. Number of uncontrollable power nodes; $P_{init}$ mode.}
\vspace{-3mm}
\label{Cluster_Control}
\end{figure}

Figure \ref{figSizeComm} shows that for a given size of power disturbance, the average yield (over 100 scenarios) decreases as the number of uncontrollable nodes increases. In addition, it shows that for the same number of uncontrollable nodes, the average yield decreases as the size of cluster increases. We observed that for the cases that all the uncontrollable areas can be stabilized using droop control; i.e. $I^k=1$, the yield of partial communication is very close to the yield of full communication. However, for the cases that at least one uncontrollable area is unstable, the yield decreases significantly. Since large clusters could contain more generators and loads, the effect of loss of such clusters is more severe. Moreover, Figure \ref{Cluster_Control} shows the impact of size and structure of uncontrollable areas on their stability. It can be seen that as the size of uncontrollable area increases the fraction of unstable cases increases. In addition, it can be seen that as the size of clusters increases, it is less probable to lose a cluster. This is due to the fact that a large cluster could contain more power, and losing it could cause a significant loss of power. Finally, observations from Figures \ref{figSizeComm} and \ref{Cluster_Control} show that every uncontrollable power node does not fail; i.e. ``point-wise" failure model is not appropriate for power-communication interdependency. 

\subsection{Effect of Size of Power Loss}
\begin{figure}[h]
\centering
\includegraphics[scale=0.055]{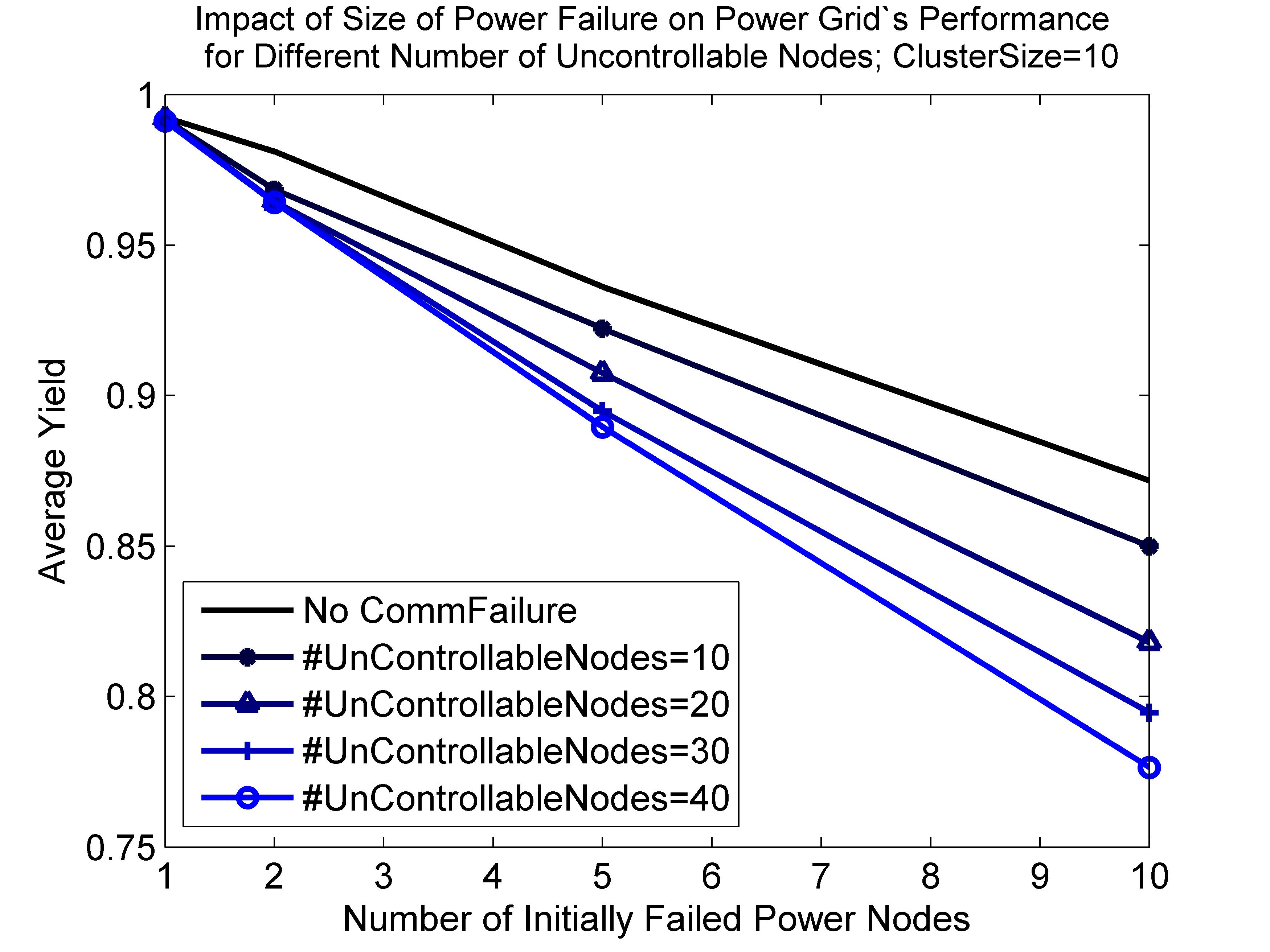}
\caption{Average yield vs. Number of initially failed power nodes; $P_{init}$ mode.}
\vspace{-3mm}
\label{figSizePow_ratio}
\end{figure}

Next, we consider the effect of size of power failure on performance. Figure \ref{figSizePow_ratio} shows that as the size of power failure increases, the average yield decreases. Moreover, it shows that the impact of communication loss increases as the number of power failures increases. This is due to the fact that for larger sizes of uncontrollable areas, it is harder to control a power disturbance and the impact on the yield is more severe.

\subsection{Effect of Different Modes}

\begin{figure}[h]
\centering
\includegraphics[scale=0.055]{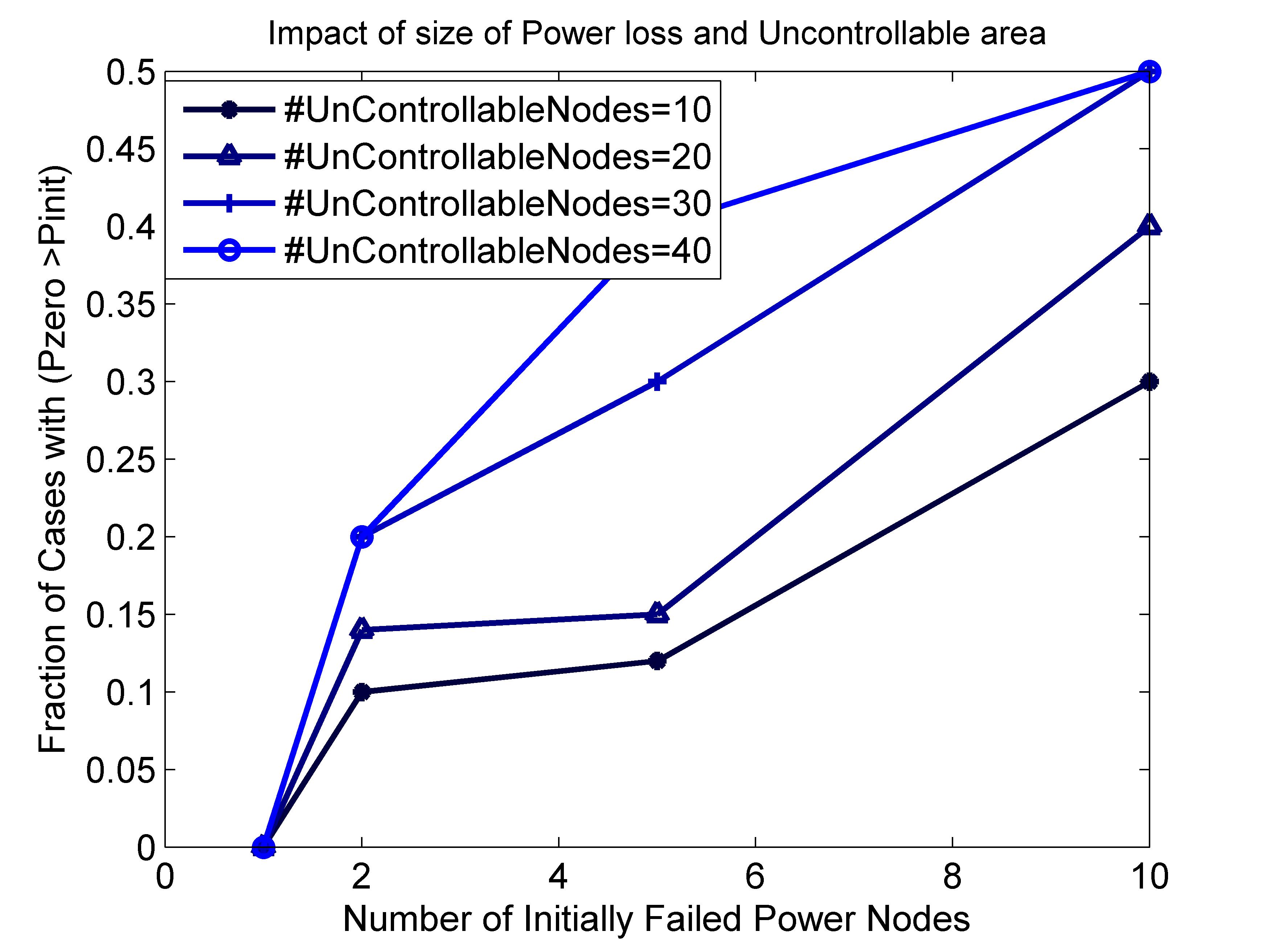}
\caption{Fraction of cases with ($Y(P_{zero})>Y(P_{init})$) vs. Number of initially failed power nodes.}
\vspace{-3mm}
\label{Mode_Ratio}
\end{figure}

Finally, we compare the impact of different predefined modes $P_{init}$ and $P_{zero}$ described in Section \ref{Partial_Emergency}. By simulating different failure scenarios under both modes, we observed that there exist scenarios where controlling the grid under $P_{zero}$ mode leads to higher yield than $P_{init}$ mode. In fact, in such scenarios, keeping the uncontrollable area operating at the $P_{init}$ mode was infeasible; thus, they were islanded due to instability. But, under $P_{zero}$ mode, it is possible to keep an uncontrollable area and use the buses in that area for transmitting power. Figure \ref{Mode_Ratio} shows that although in most cases the $P_{init}$ mode results in a higher yield, the fraction of scenarios with $Y(P_{zero})>Y(P_{init})$ increases as the number of uncontrollable nodes or failed power nodes increases. In particular, we observed that for the cases where $P_{init}$ has the higher yield, the average difference in yieldis $7\%$ and for the cases where $P_{zero}$ has the higher yield, the average difference is $8\%$.

\section{Conclusion}\label{Conclusion_Sec}
In this paper, we showed that although controlling the power grid using the communication network could be very beneficial, it could be harmful if we lose part of the control and communication network. Therefore, it is essential to have a thorough analysis on the impact of communication network on power grid to identify the vulnerabilities of the system. In particular, we showed that this impact is a function of several parameters including the size and structure of the communication loss. Therefore, it is very important to not only design a robust communication network, but also allocate the communication nodes to the power grid so that the negative impact of communication loss is minimized. Another direction of research is to design more intelligent local controllers so that in the lack of communication, the nodes can stabilize the grid even during large disturbances.

\bibliographystyle{IEEEtran}
\bibliography{reference}

\end{document}